\documentclass[12pt]{amsart}
\makeatletter

\@addtoreset{equation}{section}
\makeatother
\usepackage[margin=3cm]{geometry}
\newenvironment{nouppercase}{%
  \renewcommand{\uppercasenonmath}[1]{}}{}
\usepackage{amsmath,amssymb,amsthm,color}
\usepackage{mathrsfs} 
\usepackage{ytableau}
\usepackage[all]{xy}
\allowdisplaybreaks[1]
\theoremstyle{definition}
\newtheorem{thm}[equation]{Theorem}
 
\newtheorem{defi}[equation]{Definition} 
\newtheorem*{conj*}{Conjecture}
\newtheorem*{theorem*}{Theorem}
\newtheorem{remark}[equation]{Remark}

\newtheorem{prob}[equation]{Problem}

\newtheorem{exam}[equation]{Example}

\begin{document}

\title{Shuffle product formula of the Schur multiple zeta values of hook type}
\author{Maki Nakasuji}
\address{Department of Information and Communication Science, Faculty of Science, Sophia
University, 7-1 Kioi-cho, Chiyoda-ku, Tokyo 102-8554, Japan}
\email{ nakasuji@sophia.ac.jp}
\author[W. Takeda]{Wataru Takeda}
\address{Department of Applied Mathematics, Tokyo University of Science,
1-3 Kagurazaka, Shinjuku-ku, Tokyo 162-8601, Japan.}
\email{w.takeda@rs.tus.ac.jp}

\begin{nouppercase}
\maketitle
\end{nouppercase}
\begin{abstract}
We discuss the shuffle product of the Schur multiple zeta values, which are the special values of Schur multiple zeta functions.
We first define $2$-labeled Schur posets to generalize Yamamoto's integral expression of the multiple zeta values and consider the product of hook-type Schur multiple zeta values by using these posets. Then, for the derived terms,
we introduce a modified Hurwitz-type Schur multiple zeta function of hook type, named an elementary factorial Schur multiple zeta function. 
Furthermore, we generalize $2$-labeled Schur posets to consider the shuffle product of the elementary factorial Schur multiple zeta values and obtain 
an explicit formula for their shuffle product.
\end{abstract}
\section{Introduction}
For positive integers $r$, $k_1, k_2, \ldots, k_r$ with $k_r>1$,  
a multiple zeta value of Euler-Zagier type
is defined by
\[
 \zeta(k_1,\ldots,k_r)=\sum_{1\le n_1<\cdots< n_r}\frac{1}{n_1^{k_1}\cdots n_r^{k_r}},
 \]
 where the summation runs over all the size $r$ sets of ordered positive integers.
 This value also has an iterated integral representation:
\begin{equation}
\label{iterated}
    \zeta(k_1,\ldots,k_r)=\underset{1>x_1>\cdots>x_k>0}{\int}\prod_{i=1}^k\omega_i(x_i),
\end{equation}
where $k = k_1+k_2+\cdots+k_r$ and $\omega_i(t) = \frac{dt}{1-t}$ if $i \in \{k_r, k_{r-1}+k_r,\ldots,k_1+k_2+\cdots+k_r\}$ and $\omega_i(t) = \frac{dt}{t}$ otherwise.
From this representation, the product of multiple zeta values is deduced as follows.
By splitting the domain of integration in the appropriate way,
we can express a product of two multiple zeta values as 
\[\left(\underset{1>x_1>\cdots>x_k>0}{\int}\prod_{i=1}^k\omega_i(x_i)\right)\left(\underset{1>y_1>\cdots>y_l>0}{\int}\prod_{i=1}^l\omega_i(y_i)\right)=\sum_{(z_1,\ldots,z_{k+l})}\underset{1>z_1>\cdots>z_{k+l}>0}{\int}\prod_{i=1}^{k+l}\omega_i(z_i),
\]
where the sum $\displaystyle \sum_{(z_1,\ldots,z_{k+l})}$ means the summation of the $(k+l)$-tuples $(z_1,\ldots,z_{k+l})$ satisfying $\{z_1,\ldots,z_{k+l}\}=\{x_1,\ldots,x_k,y_1,\ldots,y_l\}$ and the partial order 
\[\alpha< \beta \text{ in } Z \Longleftrightarrow \alpha, \beta\in X \text{ and }\alpha< \beta \text{ in } X, \text{ or
} \alpha, \beta\in Y \text{ and }\alpha< \beta \text{ in }Y,\]
with $X=\{x_1,\ldots,x_{k}\},Y=\{y_1,\ldots,y_{l}\}$ and $Z=\{z_1,\ldots,z_{k+l}\}$.
Such a product is called a {\it shuffle product} and was first introduced by Eilenberg and Mac Lane in \cite[Section 5]{em}.
Since its introduction, this kind of product has been studied in the theory of multiple zeta values along with another product, called the {\it harmonic product} \cite{hoffman,hoffmanohno}.
In  \cite{y}, Yamamoto developed a graphical representation of iterated integrals in terms of labeled posets.
This representation gives a new expression of an integral representation of the multiple zeta values and
leads to a combinatorial explanation of shuffle products.
On the other hand, the shuffle product of 
the multiple zeta-star values defined by
\[
 \zeta^{\star}(k_1,\ldots,k_r)=\sum_{1\le n_1\leq \cdots\leq n_r}\frac{1}{n_1^{k_1}\cdots n_r^{k_r}}
 \]
had not been obtained explicitly, although they share Yamamoto's integral expression.
Keeping these facts in mind, in this article, we will discuss the shuffle product of the Schur multiple zeta values, which are the special values of Schur multiple zeta functions regarded as a generalization of both multiple zeta and multiple zeta-star functions of Euler-Zagier type.

The {\it Schur multiple zeta functions} introduced by Nakasuji, Phuksuwan, and Yamasaki  \cite{npy} 
are defined as sums like the usual Schur polynomials of multiple zeta functions.
More precisely,
for any partition $\lambda$, i.e., a non-increasing sequence $( \lambda_1, \ldots, \lambda_n)$ of { positive} integers,
we associate Young diagram
$D(\lambda)=\{(i, j)\in {\mathbb Z}^2 ~|~ 1\leq i\leq n, 1\leq j\leq \lambda_i\}$ depicted as a collection of square boxes with the $i$-th row having $\lambda_i$ boxes, as with matrices.
For a partition $\lambda$, a Young tableau $T=(t_{ij})$ of shape $\lambda$ over a set $X$ is a filling of $D(\lambda)$ obtained by putting $t_{ij}\in X$ into the $(i,j)$ box of $D(\lambda)$.
 We denote by $T(\lambda,X)$ the set of all Young tableaux of shape $\lambda$ over $X$ and denote by $SSYT(\lambda)$ the set of semi-standard Young tableaux {$(t_{ij})\in T(\lambda,\mathbb N)$ which satisfies} weakly increasing from left to right in each row $i$
and strictly increasing from top to bottom in each column $j$.
Then for a given set ${ \pmb s}=(s_{ij})\in T(\lambda,\mathbb{C})$ of variables, 
{\it Schur multiple zeta functions} of {shape} $\lambda$ are defined to be
\[
\zeta_{\lambda}({ \pmb s})=\sum_{M\in SSYT(\lambda)}\frac{1}{M^{ \pmb s}},
\]
where $M^{ \pmb s}=\displaystyle{\prod_{(i, j)\in D(\lambda)}m_{ij}^{s_{ij}}}$ for $M=(m_{ij})\in SSYT(\lambda)$. 
The function $\zeta_\lambda(\pmb s)$ converges absolutely in \[
 W_{\lambda}
=
\left\{{\pmb s}=(s_{ij})\in T(\lambda,\mathbb{C})\,\left|\,
\begin{array}{l}
 \text{$\Re(s_{ij})\ge 1$ for all $(i,j)\in D(\lambda) \setminus C(\lambda)$} \\[3pt]
 \text{$\Re(s_{ij})>1$ for all $(i,j)\in C(\lambda)$}
\end{array}
\right.
\right\}
\]
with $C(\lambda)$ being the set of all corners of $\lambda$. Throughout this article, we assume that all variables of $\zeta_\lambda$ are elements of $W_\lambda$. We sometimes write $\zeta_{\lambda}({ \pmb s})$ as $\pmb s$ for short if there is no confusion.

In formulating the shuffle product of the Schur multiple zeta values (SMZVs), a similar difficulty as in the product of the multiple zeta-star values occur, even if the shape of Schur multiple zeta values is restricted to a simple one, say the hook type, which is $\lambda=(\lambda_1, \{1\}^r)$ ($\lambda_1\geq 1$), where $\{1\}^r$ means $\underbrace{1,\ldots,1}_r$. 
However, to overcome this difficulty, we define a new generalization of Schur multiple zeta functions of hook type which is modified Hurwitz-type Schur multiple zeta functions, and are able to give an expression for the product. Moreover, we obtain an explicit form of the shuffle product of these functions.

This article is organized as follows. In Section \ref{preex}, we will show the basic idea with some examples, to introduce the shuffle product of Schur multiple zeta values of hook type.
In Section \ref{2posets}, we will define $2$-labeled Schur posets to generalize Yamamoto's integral expression of multiple zeta values and consider the product of  Schur multiple zeta values of hook type by using these posets.
In Section \ref{elementarysec}, we will introduce the elementary factorial Schur multiple zeta functions, which are modified Hurwitz-type Schur multiple zeta functions of hook type, and rewrite the examples in Sections \ref{preex} and \ref{2posets}. Then we will generalize $2$-labeled Schur posets to consider the shuffle product of the elementary factorial Schur multiple zeta values and obtain 
an explicit form of their shuffle product (Theorem \ref{shuffle3}).
In the final section, we will present further problems and give brief examples.

\section{Preliminary examples}
\label{preex}
In this section, we will lay out some basic ideas in order to introduce the shuffle product of the hook-type Schur multiple zeta functions through examples.
Throughout this paper, we use notation $e^r=(\{1\}^r)$ and $h^r=(r)$.
\begin{exam}
\label{firstex}
Let $\lambda=(2,1)$. Then it holds that
\begin{align*}
\ytableausetup{boxsize=normal,aligntableaux=center}
&\begin{ytableau}
  1&3\\
  2\\
\end{ytableau}\times\begin{ytableau}
  2
\end{ytableau}\\
&=\sum_{\substack{(m_{ij})\in SSYT(\lambda)\\ (m)\in SSYT(e^1)}}\left[\frac{1}{m^2(m+m_{11})m_{12}^3(m+m_{21})^2}+\frac{1}{m(m+m_{11})^2m_{12}^3(m+m_{21})^2}\right.\\
&\left.+\frac{2}{m(m+m_{11})m_{12}^3(m+m_{21})^3}\right]+2\ \begin{ytableau}
  1&3\\
  2\\
  2
\end{ytableau}+4\ \begin{ytableau}
  1&3\\
  1\\
  3
\end{ytableau}.
\end{align*}  

\end{exam}
\begin{proof}
From \cite{nn} and \cite{npy}, we find that 
\[\begin{ytableau}
  1&3\\
  2\\
\end{ytableau}=\int_{0}^1\frac{dz_1}{z_1}\int_{0}^{z_1}\frac{dz_2}{1-z_2}\int_{0}^{z_2}\frac{dz_3}{1-z_3}\sum_{m_{12}=1}^\infty\frac{1-z_3^{m_{12}}}{m_{12}^{3}}.\]
Also, it is well known that 
\[\begin{ytableau}
 2
\end{ytableau}=\int_{0}^1\frac{dw_1}{w_1}\int_{0}^{w_1}\frac{dw_2}{1-w_2}.\]
Considering the order of $z_1,z_2,z_3,w_1,w_2$, we find that
\begin{align*}
    \begin{ytableau}
  1&3\\
  2\\
\end{ytableau}\times\begin{ytableau}
  2
\end{ytableau}=&\int_{0}^1\frac{dz_1}{z_1}\int_{0}^{z_1}\frac{dz_2}{1-z_2}\int_{0}^{z_2}\frac{dz_3}{1-z_3}\int_{0}^{z_3}\frac{dw_1}{w_1}\int_{0}^{w_1}\frac{dw_2}{1-w_2}\sum_{m_{12}=1}^\infty\frac{1-z_3^{m_{12}}}{m_{12}^{3}}\\
&+\int_{0}^1\frac{dz_1}{z_1}\int_{0}^{z_1}\frac{dz_2}{1-z_2}\int_{0}^{z_2}\frac{dw_1}{w_1}\int_{0}^{w_1}\frac{dz_3}{1-z_3}\int_{0}^{z_3}\frac{dw_2}{1-w_2}\sum_{m_{12}=1}^\infty\frac{1-z_3^{m_{12}}}{m_{12}^{3}}\\
&+\int_{0}^1\frac{dz_1}{z_1}\int_{0}^{z_1}\frac{dw_1}{w_1}\int_{0}^{w_1}\frac{dz_2}{1-z_2}\int_{0}^{z_2}\frac{dz_3}{1-z_3}\int_{0}^{z_3}\frac{dw_2}{1-w_2}\sum_{m_{12}=1}^\infty\frac{1-z_3^{m_{12}}}{m_{12}^{3}}\\
&+\int_{0}^1\frac{dw_1}{w_1}\int_{0}^{w_1}\frac{dz_1}{z_1}\int_{0}^{z_1}\frac{dz_2}{1-z_2}\int_{0}^{z_2}\frac{dz_3}{1-z_3}\int_{0}^{z_3}\frac{dw_2}{1-w_2}\sum_{m_{12}=1}^\infty\frac{1-z_3^{m_{12}}}{m_{12}^{3}}\\
&+\int_{0}^1\frac{dz_1}{z_1}\int_{0}^{z_1}\frac{dz_2}{1-z_2}\int_{0}^{z_2}\frac{dw_1}{w_1}\int_{0}^{w_1}\frac{dw_2}{1-w_2}\int_{0}^{w_2}\frac{dz_3}{1-z_3}\sum_{m_{12}=1}^\infty\frac{1-z_3^{m_{12}}}{m_{12}^{3}}\\
&+\int_{0}^1\frac{dz_1}{z_1}\int_{0}^{z_1}\frac{dw_1}{w_1}\int_{0}^{w_1}\frac{dz_2}{1-z_2}\int_{0}^{z_2}\frac{dw_2}{1-w_2}\int_{0}^{w_2}\frac{dz_3}{1-z_3}\sum_{m_{12}=1}^\infty\frac{1-z_3^{m_{12}}}{m_{12}^{3}}\\
&+\int_{0}^1\frac{dw_1}{w_1}\int_{0}^{w_1}\frac{dz_1}{z_1}\int_{0}^{z_1}\frac{dz_2}{1-z_2}\int_{0}^{z_2}\frac{dw_2}{1-w_2}\int_{0}^{w_2}\frac{dz_3}{1-z_3}\sum_{m_{12}=1}^\infty\frac{1-z_3^{m_{12}}}{m_{12}^{3}}\\
&+\int_{0}^1\frac{dz_1}{z_1}\int_{0}^{z_1}\frac{dw_1}{w_1}\int_{0}^{w_1}\frac{dw_2}{1-w_2}\int_{0}^{w_2}\frac{dz_2}{1-z_2}\int_{0}^{z_2}\frac{dz_3}{1-z_3}\sum_{m_{12}=1}^\infty\frac{1-z_3^{m_{12}}}{m_{12}^{3}}\\
&+\int_{0}^1\frac{dw_1}{w_1}\int_{0}^{w_1}\frac{dz_1}{z_1}\int_{0}^{z_1}\frac{dw_2}{1-w_2}\int_{0}^{w_2}\frac{dz_2}{1-z_2}\int_{0}^{z_2}\frac{dz_3}{1-z_3}\sum_{m_{12}=1}^\infty\frac{1-z_3^{m_{12}}}{m_{12}^{3}}\\
&+\int_{0}^1\frac{dw_1}{w_1}\int_{0}^{w_1}\frac{dw_2}{1-w_2}\int_{0}^{w_2}\frac{dz_1}{z_1}\int_{0}^{z_1}\frac{dz_2}{1-z_2}\int_{0}^{z_2}\frac{dz_3}{1-z_3}\sum_{m_{12}=1}^\infty\frac{1-z_3^{m_{12}}}{m_{12}^{3}}\\
&=:\sum_{k=1}^{10}I_k.
\end{align*}
We note that for $|z|<1$ we have $\int_{0}^{t}\frac{dz}{1-z}-\log(1-t)=\sum_{m=1}^{\infty}\frac{t^m}{m}$.
By computing directly, we obtain
\begin{align*}
    I_1&=\int_{0}^1\frac{dz_1}{z_1}\int_{0}^{z_1}\frac{dz_2}{1-z_2}\int_{0}^{z_2}\frac{dz_3}{1-z_3}\int_{0}^{z_3}\frac{dw_1}{w_1}\int_{0}^{w_1}\frac{dw_2}{1-w_2}\sum_{m_{12}=1}^\infty\frac{1-z_3^{m_{12}}}{m_{12}^{3}}\\
    &=\int_{0}^1\frac{dz_1}{z_1}\int_{0}^{z_1}\frac{dz_2}{1-z_2}\int_{0}^{z_2}\frac{dz_3}{1-z_3}\int_{0}^{z_3}\frac{dw_1}{w_1}\sum_{m=1}^{\infty}\frac{w_1^m}{m}\sum_{m_{12}=1}^\infty\frac{1-z_3^{m_{12}}}{m_{12}^{3}}\\
    &=\int_{0}^1\frac{dz_1}{z_1}\int_{0}^{z_1}\frac{dz_2}{1-z_2}\int_{0}^{z_2}\frac{dz_3}{1-z_3}\sum_{m=1}^{\infty}\frac{z_3^m}{m^2}\sum_{m_{12}=1}^\infty\frac{1-z_3^{m_{12}}}{m_{12}^{3}}.
    \intertext{Since both series converge absolutely in $|z_3|<1$, we find that }
    I_1&=\int_{0}^1\frac{dz_1}{z_1}\int_{0}^{z_1}\frac{dz_2}{1-z_2}\int_{0}^{z_2}\frac{dz_3}{1-z_3}\sum_{m=1}^{\infty}\sum_{m_{12}=1}^\infty\frac{z_3^m(1-z_3^{m_{12}})}{m^2m_{12}^{3}}\\
    &=\int_{0}^1\frac{dz_1}{z_1}\int_{0}^{z_1}\frac{dz_2}{1-z_2}\sum_{m_{11}=1}^{m_{12}}\sum_{m=1}^{\infty}\sum_{m_{12}=1}^\infty\frac{z_2^{m+m_{11}}}{m^2(m+m_{11})m_{12}^{3}}\\
   &=\int_{0}^1\frac{dz_1}{z_1}\sum_{m_{21}=m_{11}+1}^{\infty}\sum_{m_{11}=1}^{m_{12}}\sum_{m=1}^{\infty}\sum_{m_{12}=1}^\infty\frac{z_1^{m+m_{21}}}{m^2(m+m_{11})(m+m_{21})m_{12}^{3}}\\
   &=\sum_{m_{21}=m_{11}+1}^{\infty}\sum_{m_{11}=1}^{m_{12}}\sum_{m=1}^{\infty}\sum_{m_{12}=1}^\infty\frac{1}{m^2(m+m_{11})(m+m_{21})^2m_{12}^{3}}\\
   &=\sum_{\substack{(m_{ij})\in SSYT(\lambda)\\ (m)\in SSYT(e^1)}}\frac{1}{m^2(m+m_{11})(m+m_{21})^2m_{12}^{3}}.
\end{align*}
Similarly, we can compute $I_2,I_3,I_4$ as above and obtain
\begin{align*}
    I_2&=\sum_{\substack{(m_{ij})\in SSYT(\lambda)\\ (m)\in SSYT(e^1)}}\frac{1}{m(m+m_{11})^2(m+m_{21})^2m_{12}^{3}},\\
    I_3=I_4&=\sum_{\substack{(m_{ij})\in SSYT(\lambda)\\ (m)\in SSYT(e^1)}}\frac{1}{m(m+m_{11})(m+m_{21})^3m_{12}^{3}}.
\end{align*}
On the other hand, we find that 
\begin{align*}
    I_5=I_{10}&=\sum_{\substack{(m_{ij})\in SSYT(\lambda)\\ (m)\in SSYT(e^1)}}\frac{1}{m_{11}(m+m_{11})^2(m+m_{21})^2m_{12}^{3}}=\begin{ytableau}
  1&3\\
  2\\
  2
\end{ytableau},\\
\intertext{and}
    I_6=\cdots=I_9&=\sum_{\substack{(m_{ij})\in SSYT(\lambda)\\ (m)\in SSYT(e^1)}}\frac{1}{m_{11}(m+m_{11})(m+m_{21})^3m_{12}^{3}}=\begin{ytableau}
  1&3\\
  1\\
  3
\end{ytableau}.
\end{align*}
This completes the proof of this example.
\end{proof}
The above example shows that non-Schur multiple zeta values appear in the product of a Schur multiple zeta value of hook type and a multiple zeta value of Euler-Zagier type. We will consider non-Schur multiple zeta values in Section \ref{elementarysec} and confirm that they can be expressed as sums of Schur multiple zeta values. In this section, we leave them in the form of infinite sums without any particular definition.
Focusing on the proof of Example \ref{firstex}, we can apply the same calculation for the shuffle product $\zeta_{(r,1)}(\pmb k)$ and a multiple zeta value of Euler-Zagier type, and obtain the following example.
\begin{exam}
\label{12112}
Let $\lambda=(r,1)$ and
$\begin{ytableau}
  \pmb k
\end{ytableau}=\begin{ytableau}
  k_2&\cdots& k_r
\end{ytableau}$ with $k_r>1$. Then it holds that
\begin{align*}\ytableausetup{boxsize=normal,aligntableaux=center}
&\begin{ytableau}
  1&\pmb k\\
  2
\end{ytableau}\times\begin{ytableau}
  1\\
  1\\
  2
\end{ytableau}\\
&=\sum_{\substack{(m_{ij})\in SSYT(\lambda)\\ (m_i)\in SSYT(e^3)}}\prod_{j=2}^r\frac{1}{m_{1j}^{k_{j}}}\left[\frac{1}{m_1m_2m_3^2(m_3+m_{11})(m_3+m_{21})^2}+\frac{1}{m_1m_2m_3(m_3+m_{11})^2(m_3+m_{21})^2}\right.\\
&\left.+\frac{2}{m_1m_2m_3(m_3+m_{11})(m_3+m_{21})^3}+\frac{2}{m_1m_2(m_2+m_{11})(m_3+m_{11})^2(m_3+m_{21})^2}\right.\\
&\left.+\frac{4}{m_1m_2(m_2+m_{11})(m_3+m_{11})(m_3+m_{21})^3}+\frac{3}{m_1(m_1+m_{11})(m_2+m_{11})(m_3+m_{11})^2(m_3+m_{21})^2}\right.\\
&\left.+\frac{6}{m_1(m_1+m_{11})(m_2+m_{11})(m_3+m_{11})(m_3+m_{21})^3}\right.\\
&\left.+\frac{1}{m_1(m_1+m_{11})(m_2+m_{11})^2(m_3+m_{11})(m_3+m_{21})^2}\right]\\
&+\begin{ytableau}
  1&\pmb k\\
  2\\
  1\\
  1\\
  2
\end{ytableau}+2\ \begin{ytableau}
  1&\pmb k\\
  1\\
  2\\
  1\\
  2
\end{ytableau}+4\ \begin{ytableau}
  1&\pmb k\\
  1\\
  1\\
  2\\
  2
\end{ytableau}+8\ \begin{ytableau}
  1&\pmb k\\
  1\\
  1\\
  1\\
  3
\end{ytableau}.
\end{align*}  
\end{exam}
\section{2-posets and associated integral}
\label{2posets}
In \cite{y}, Yamamoto gave a graphical expression of an iterated integral representation of the multiple zeta values of Euler-Zagier type. This has been used in various contexts in the theory of multiple zeta values and has led to many remarkable combinatorial proofs. In this section, we introduce {\it $2$-labeled Schur posets} to generalize Yamamoto's expression as follows.
\begin{defi}
\label{defi:sp}
 Let $\begin{ytableau}
  \pmb k
\end{ytableau}=\begin{ytableau}
  k_2&\cdots& k_r
\end{ytableau}$.
\begin{enumerate}
\item A $2$-labeled Schur poset is a pair $(X,\pmb k)=(X,\pmb k;\delta_X)$ consisting of 
a poset $X$, a Young tableau $\pmb k$, and a map $\delta_X\colon X\to\{0,1\}$ called a {labeling map}. 
\item A $2$-labeled Schur poset $(X,\pmb k)$ is said to be {admissible} 
if $\delta_X(x)=1$ for all minimal $x\in X$ and 
$\delta_X(x)=0$ for all maximal $x\in X$. 
\item For any poset $X$, we put 
\[\Delta(X):=\bigl\{(t_x)_{x\in X}\in[0,1]^X\bigm| 
t_x<t_y \text{ if } x<y\bigr\}. \]
\item For admissible $2$-labeled Schur posets $I(X,\pmb k)$, 
we define the associated integral by 
\begin{equation}\label{eq:I}
I(X,\pmb k):=\int_{\Delta(X)}\prod_{x\in X}\omega_{\delta_X(x)}(t_x)\sum_{m_{1j}\in SSYT(h_{r-1})}\frac{1-t_{x}^{m_{12}}}{m_{12}^{k_{2}}\cdots m_{1r}^{k_r}},
\end{equation}
where $x$ is the minimal $x\in X$.
Here $\omega_0(t)=\frac{dt}{t}$ and $\omega_1(t)=\frac{dt}{1-t}$.
\end{enumerate}
\end{defi}
 We use Hasse diagrams to indicate $2$-labeled Schur posets, 
with vertices $\circ$ and $\bullet$ corresponding to $\delta(x)=0$ and $1$, 
respectively. For example, $X=\{x_1,x_2,x_3,x_4\}$ with order $x_1>x_2>x_3>x_4$ and label $(\delta(x_1),\delta(x_2),\delta(x_3),\delta(x_4))=(0,1,0,1)$ is represented as the diagram
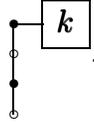
\[\begin{xy}
(7,4)*{\begin{ytableau}
  \pmb k
\end{ytableau}},
{(0,4) \ar @{-} (4,4)},
{(0,4) \ar @{{*}-o} (0,0)},
{(0,0) \ar @{-{*}} (0,-4)},
{(0,-4) \ar @{-o} (0,-8)},
\end{xy}.\]
This diagram corresponds to
\[
\begin{ytableau}
  2&\pmb k\\
  2
\end{ytableau}=\int_{1>x_1>x_2>x_3>x_4>0}\frac{dx_1}{x_1}\frac{dx_2}{1-x_2}\frac{dx_3}{x_3}\frac{dx_4}{1-x_4}\sum_{m_{1j}\in SSYT(h_{r-1})}\frac{1-x_4^{m_{12}}}{m_{12}^{k_{2}}\cdots m_{1r}^{k_r}}.
\]
    If $\pmb k=\emptyset$, we assign $1$ to the inner empty
sum. Then admissible $2$-labeled Schur posets  $(X,\emptyset)$ correspond to multiple zeta values of Euler-Zagier type. For example, when $X=\{x_1,x_2,x_3,x_4,x_5\}$ with order $x_1>x_2>x_3>x_4>x_5$ and label $(\delta(x_1),\delta(x_2),\delta(x_3),\delta(x_4),\delta(x_5))=(0,1,0,0,1)$, the $2$-labeled Schur poset $(X,\emptyset;\delta_X)$ is represented as diagram
\[\begin{xy}
(7,4)*{\emptyset},
{(0,4) \ar @{-} (4,4)},
{(0,4) \ar @{{*}-o} (0,0)},
{(0,0) \ar @{-o} (0,-4)},
{(0,-4) \ar @{-{*}} (0,-8)},
{(0,-8) \ar @{-o} (0,-12)},
\end{xy}\longleftrightarrow\begin{xy}
{(0,4) \ar @{{*}-o} (0,0)},
{(0,0) \ar @{-o} (0,-4)},
{(0,-4) \ar @{-{*}} (0,-8)},
{(0,-8) \ar @{-o} (0,-12)},
\end{xy}\ \  \text{ and }I(X,\emptyset)=\int_{x_1>x_2>x_3>x_4>x_5}\frac{dx_1}{x_1}\frac{dx_2}{1-x_2}\frac{dx_3}{x_3}\frac{dx_4}{x_4}\frac{dx_5}{1-x_5}=\zeta(3,2).\]
Here, we note that the way the Hasse diagram is expressed is opposite to those of Yamamoto \cite{y} and other papers using Yamamoto's integral expression. One may think that the opposite order is natural to compare Schur multiple zeta functions to Young diagrams.
For $2$-labeled Schur posets $(X,\pmb k)$ and $(Y,\emptyset)$  one can define their direct sum $(X\amalg Y,\pmb k)$: Its
underlying poset is the direct sum of finite sets $X$ and $Y$ endowed with the partial order 
\[x< y \text{ in } X\amalg Y \Leftrightarrow x, y\in X \text{ and }x< y \text{ in } X \text{ or
} x, y\in Y \text{ and }x< y \text{ in }Y.\]
\begin{thm}[Shuffle product of SMZV of hook type]
\label{shuffle}
If $2$-labeled Schur posets $(X,\pmb k)$ and $(Y,\emptyset)$ are admissible, then their direct sum $X\amalg Y$ is also admissible and
\[I(X\amalg Y,\pmb k)=I(X,\pmb k)I(Y,\emptyset).\]
\end{thm}
When we use a Hasse diagram to compute a product \[I\left(\ \begin{xy}
(-3,4)*{x_4},
(-3,0)*{x_3},
(-3,-4)*{x_2},
(-3,-8)*{x_1},
{(0,4) \ar @{{*}-o} (0,0)},
{(0,0) \ar @{-{*}} (0,-4)},
{(0,-4) \ar @{-o} (0,-8)},
\end{xy}\ ,\ \pmb k\right)I\left(\ \begin{xy}
(-3,2)*{y_2},
(-3,-2)*{y_1},
{(0,2) \ar @{{*}-o} (0,-2)}
\end{xy}\ ,\ \emptyset\right),\] we find that
\begin{align*}
    I\left(\begin{xy}
(-3,10)*{y_2},
(-3,6)*{y_1},
(-3,2)*{x_4},
(-3,-2)*{x_3},
(-3,-6)*{x_2},
(-3,-10)*{x_1},
{(0,10) \ar @{{*}-o} (0,6)},
{(0,6) \ar @{-{*}} (0,2)},
{(0,2) \ar @{-o} (0,-2)},
{(0,-2) \ar @{-{*}} (0,-6)},
{(0,-6) \ar @{-o} (0,-10)},
\end{xy}\ ,\pmb k\right)&=\sum_{\substack{(m_{ij})\in SSYT(\lambda)\\ (m)\in SSYT(e^1)}}\prod_{j=2}^r\frac{1}{m_{1j}^{k_{j}}}\frac{1}{m^2(m+m_{11})^2(m+m_{21})^2},\\
I\left(\begin{xy}
(-3,10)*{x_4},
(-3,6)*{x_3},
(-3,2)*{x_2},
(-3,-2)*{x_1},
(-3,-6)*{y_2},
(-3,-10)*{y_1},
{(0,10) \ar @{{*}-o} (0,6)},
{(0,6) \ar @{-{*}} (0,2)},
{(0,2) \ar @{-o} (0,-2)},
{(0,-2) \ar @{-{*}} (0,-6)},
{(0,-6) \ar @{-o} (0,-10)},
\end{xy}\ ,\pmb k\right)&=\begin{ytableau}
  2&\pmb k\\
  2\\
  2
\end{ytableau}.
\end{align*}
Therefore, we have to distinguish the poset $(X,\pmb k)$ corresponding to the Schur multiple zeta value from $(Y,\emptyset)$, which is associated with the multiple zeta values of Euler-Zagier type.
For this purpose, we use vertices $\{\circ,\bullet\}$ and $\{\square,\blacksquare\}$ for $(X,\pmb k)$ and $(Y,\emptyset)$, respectively. For example, 
\begin{align*}
   \begin{xy}
(7,4)*{\begin{ytableau}
  \pmb k
\end{ytableau}},
{(0,4) \ar @{-} (4,4)},
{(0,4) \ar @{{*}-o} (0,0)},
{(0,0) \ar @{-{*}} (0,-4)},
{(0,-4) \ar @{-o} (0,-8)},
\end{xy}&\longleftrightarrow\begin{ytableau}
  2&\pmb k\\
  2
\end{ytableau}=\int_{1>x_1>x_2>x_3>x_4>0}\frac{dx_1}{x_1}\frac{dx_2}{1-x_2}\frac{dx_3}{x_3}\frac{dx_4}{1-x_4}\sum_{m_{1j}\in SSYT(h_{r-1})}\frac{1-x_4^{m_{12}}}{m_{12}^{k_{2}}\cdots m_{1r}^{k_r}},\\
\begin{xy}
(0,5)*++{\scriptstyle \blacksquare},
{(0,4) \ar @{-} (0,0)},
(0,0)*++{\scriptstyle \square},
{(0,0) \ar @{-} (0,-4)},
(0,-4)*++{\scriptstyle \square},
{(0,-4) \ar @{-} (0,-8)},
(0,-8)*++{\scriptstyle \square},
\end{xy}&\longleftrightarrow\begin{ytableau}
 4
\end{ytableau}=\int_{1>y_1>y_2>y_3>y_4>0}\frac{dy_1}{y_1}\frac{dy_2}{y_2}\frac{dy_3}{y_3}\frac{dy_4}{1-y_4}.
\end{align*}
\begin{exam}
\label{222}
\begin{align*}
&\begin{ytableau}
  2&2\\
  2\\
\end{ytableau}\cdot\begin{ytableau}
  2
\end{ytableau}\\
&=\sum_{\substack{(m_{ij})\in SSYT(\lambda)\\ (m)\in SSYT(e^1)}}\left[\frac{1}{m^2(m+m_{11})^2(m+m_{21})^2m_{12}^2}+\frac{2}{m(m+m_{11})^3(m+m_{21})^2m_{12}^2}\right.\\
&+\left.\frac{2}{m(m+m_{11})^2(m+m_{21})^3m_{12}^2}\right]+2\ \begin{ytableau}
  2&2\\
  2\\
  2
\end{ytableau}+2\ \begin{ytableau}
  1&2\\
  3\\
  2
\end{ytableau}+2\ \begin{ytableau}
  1&2\\
  2\\
  3
\end{ytableau}+4\ \begin{ytableau}
  2&2\\
  1\\
  3
\end{ytableau}.
\end{align*}
\end{exam}
\begin{proof}
By considering the product of posets, we have
\begin{align*}
\begin{xy}
(7,4)*{\begin{ytableau}
  2
\end{ytableau}},
{(0,4) \ar @{-} (4,4)},
{(0,4) \ar @{{*}-o} (0,0)},
{(0,0) \ar @{-{*}} (0,-4)},
{(0,-4) \ar @{-o} (0,-8)},
\end{xy}\amalg\begin{xy}
(0,1)*++{\scriptstyle \blacksquare},
{(0,0) \ar @{-} (0,-4)},
(0,-4)*++{\scriptstyle \square},
\end{xy}&=\begin{xy}
(7,4)*{\begin{ytableau}
  2
\end{ytableau}},
{(0,4) \ar @{-} (4,4)},
(0,4)*++{\scriptstyle \blacksquare},
{(0,4) \ar @{-} (0,0)},
(0,0)*++{\scriptstyle \square},
{(0,0) \ar @{-{*}} (0,-4)},
{(0,-4) \ar @{-o} (0,-8)},
{(0,-8) \ar @{-{*}} (0,-12)},
{(0,-12) \ar @{-o} (0,-16)},
\end{xy}+\begin{xy}
(7,4)*{\begin{ytableau}
  2
\end{ytableau}},
{(0,4) \ar @{-} (4,4)},
(0,4)*++{\scriptstyle \blacksquare},
{(0,4) \ar @{-{*}} (0,0)},
(0,-4)*++{\scriptstyle \square},
{(0,0) \ar @{-} (0,-4)},
{(0,-4) \ar @{-o} (0,-8)},
{(0,-8) \ar @{-{*}} (0,-12)},
{(0,-12) \ar @{-o} (0,-16)},
\end{xy}+\begin{xy}
(7,4)*{\begin{ytableau}
  2
\end{ytableau}},
{(0,4) \ar @{-} (4,4)},
(0,4)*++{\scriptstyle \blacksquare},
{(0,4) \ar @{-{*}} (0,0)},
(0,-8)*++{\scriptstyle \square},
{(0,0) \ar @{-o} (0,-4)},
{(0,-4) \ar @{-} (0,-8)},
{(0,-8) \ar @{-{*}} (0,-12)},
{(0,-12) \ar @{-o} (0,-16)},
\end{xy}+\begin{xy}
(7,4)*{\begin{ytableau}
  2
\end{ytableau}},
{(0,4) \ar @{-} (4,4)},
(0,4)*++{\scriptstyle \blacksquare},
{(0,4) \ar @{-{*}} (0,0)},
(0,-12)*++{\scriptstyle \square},
{(0,0) \ar @{-o} (0,-4)},
{(0,-4) \ar @{-{*}} (0,-8)},
{(0,-8) \ar @{-} (0,-12)},
{(0,-12) \ar @{-o} (0,-16)},
\end{xy}+\begin{xy}
(7,4)*{\begin{ytableau}
  2
\end{ytableau}},
{(0,4) \ar @{-} (4,4)},
(0,4)*++{\scriptstyle \blacksquare},
{(0,4) \ar @{-{*}} (0,0)},
(0,-16)*++{\scriptstyle \square},
{(0,0) \ar @{-o} (0,-4)},
{(0,-4) \ar @{-{*}} (0,-8)},
{(0,-8) \ar @{-o} (0,-12)},
{(0,-12) \ar @{-} (0,-16)},
\end{xy}\\
&+\begin{xy}
(7,4)*{\begin{ytableau}
  2
\end{ytableau}},
{(0,4) \ar @{-} (4,4)},
(0,0)*++{\scriptstyle \blacksquare},
{(0,4) \ar @{{*}-} (0,0)},
(0,-4)*++{\scriptstyle \square},
{(0,0) \ar @{-} (0,-4)},
{(0,-4) \ar @{-o} (0,-8)},
{(0,-8) \ar @{-{*}} (0,-12)},
{(0,-12) \ar @{-o} (0,-16)},
\end{xy}+\begin{xy}
(7,4)*{\begin{ytableau}
  2
\end{ytableau}},
{(0,4) \ar @{-} (4,4)},
(0,0)*++{\scriptstyle \blacksquare},
{(0,4) \ar @{{*}-} (0,0)},
(0,-8)*++{\scriptstyle \square},
{(0,0) \ar @{-o} (0,-4)},
{(0,-4) \ar @{-} (0,-8)},
{(0,-8) \ar @{-{*}} (0,-12)},
{(0,-12) \ar @{-o} (0,-16)},
\end{xy}+\begin{xy}
(7,4)*{\begin{ytableau}
  2
\end{ytableau}},
{(0,4) \ar @{-} (4,4)},
(0,0)*++{\scriptstyle \blacksquare},
{(0,4) \ar @{{*}-} (0,0)},
(0,-12)*++{\scriptstyle \square},
{(0,0) \ar @{-o} (0,-4)},
{(0,-4) \ar @{-{*}} (0,-8)},
{(0,-8) \ar @{-} (0,-12)},
{(0,-12) \ar @{-o} (0,-16)},
\end{xy}+\begin{xy}
(7,4)*{\begin{ytableau}
  2
\end{ytableau}},
{(0,4) \ar @{-} (4,4)},
(0,0)*++{\scriptstyle \blacksquare},
{(0,4) \ar @{{*}-} (0,0)},
(0,-16)*++{\scriptstyle \square},
{(0,0) \ar @{-o} (0,-4)},
{(0,-4) \ar @{-{*}} (0,-8)},
{(0,-8) \ar @{-o} (0,-12)},
{(0,-12) \ar @{-} (0,-16)},
\end{xy}+\begin{xy}
(7,4)*{\begin{ytableau}
  2
\end{ytableau}},
{(0,4) \ar @{-} (4,4)},
(0,-4)*++{\scriptstyle \blacksquare},
{(0,4) \ar @{{*}-o} (0,0)},
(0,-8)*++{\scriptstyle \square},
{(0,0) \ar @{-} (0,-4)},
{(0,-4) \ar @{-} (0,-8)},
{(0,-8) \ar @{-{*}} (0,-12)},
{(0,-12) \ar @{-o} (0,-16)},
\end{xy}\\
&+\begin{xy}
(7,4)*{\begin{ytableau}
  2
\end{ytableau}},
{(0,4) \ar @{-} (4,4)},
(0,-4)*++{\scriptstyle \blacksquare},
{(0,4) \ar @{{*}-o} (0,0)},
(0,-12)*++{\scriptstyle \square},
{(0,0) \ar @{-} (0,-4)},
{(0,-4) \ar @{-{*}} (0,-8)},
{(0,-8) \ar @{-} (0,-12)},
{(0,-12) \ar @{-o} (0,-16)},
\end{xy}+\begin{xy}
(7,4)*{\begin{ytableau}
  2
\end{ytableau}},
{(0,4) \ar @{-} (4,4)},
(0,-4)*++{\scriptstyle \blacksquare},
{(0,4) \ar @{{*}-o} (0,0)},
(0,-16)*++{\scriptstyle \square},
{(0,0) \ar @{-} (0,-4)},
{(0,-4) \ar @{-{*}} (0,-8)},
{(0,-8) \ar @{-o} (0,-12)},
{(0,-12) \ar @{-} (0,-16)},
\end{xy}+\begin{xy}
(7,4)*{\begin{ytableau}
  2
\end{ytableau}},
{(0,4) \ar @{-} (4,4)},
(0,-8)*++{\scriptstyle \blacksquare},
{(0,4) \ar @{{*}-o} (0,0)},
(0,-12)*++{\scriptstyle \square},
{(0,0) \ar @{-{*}} (0,-4)},
{(0,-4) \ar @{-} (0,-8)},
{(0,-8) \ar @{-} (0,-12)},
{(0,-12) \ar @{-o} (0,-16)},
\end{xy}+\begin{xy}
(7,4)*{\begin{ytableau}
  2
\end{ytableau}},
{(0,4) \ar @{-} (4,4)},
(0,-8)*++{\scriptstyle \blacksquare},
{(0,4) \ar @{{*}-o} (0,0)},
(0,-16)*++{\scriptstyle \square},
{(0,0) \ar @{-{*}} (0,-4)},
{(0,-4) \ar @{-} (0,-8)},
{(0,-8) \ar @{-o} (0,-12)},
{(0,-12) \ar @{-} (0,-16)},
\end{xy}+\begin{xy}
(7,4)*{\begin{ytableau}
  2
\end{ytableau}},
{(0,4) \ar @{-} (4,4)},
(0,-12)*++{\scriptstyle \blacksquare},
{(0,4) \ar @{{*}-o} (0,0)},
(0,-16)*++{\scriptstyle \square},
{(0,0) \ar @{-{*}} (0,-4)},
{(0,-4) \ar @{-o} (0,-8)},
{(0,-8) \ar @{-} (0,-12)},
{(0,-12) \ar @{-} (0,-16)},
\end{xy}.\\
\intertext{This leads to}
\begin{ytableau}
  2&2\\
  2\\
\end{ytableau}\cdot\begin{ytableau}
  2
\end{ytableau}&=\sum_{\substack{(m_{ij})\in SSYT(\lambda)\\ (m)\in SSYT(e^1)}}\left[\frac{1}{m^2(m+m_{11})^2(m+m_{21})^2m_{12}^2}+\frac{2}{m(m+m_{11})^3(m+m_{21})^2m_{12}^2}\right.\\
&+\left.\frac{2}{m(m+m_{11})^2(m+m_{21})^3m_{12}^2}\right]+2\ \begin{ytableau}
  2&2\\
  2\\
  2
\end{ytableau}+2\ \begin{ytableau}
  1&2\\
  3\\
  2
\end{ytableau}+2\ \begin{ytableau}
  1&2\\
  2\\
  3
\end{ytableau}+4\ \begin{ytableau}
  2&2\\
  1\\
  3
\end{ytableau}.
\end{align*}

\end{proof}
In \cite{nt}, we obtained the Pieri formula for the hook-type Schur multiple zeta functions, which gives
\[\begin{ytableau}
  2&2\\
  2\\
\end{ytableau}\cdot\begin{ytableau}
  2
\end{ytableau}=\begin{ytableau}
  2&2&2\\
  2\\
\end{ytableau}+\begin{ytableau}
  2&2\\
  2&2
\end{ytableau}+\begin{ytableau}
  2&2\\
  2\\
  2
\end{ytableau}.\]
Combining with this identity with Example \ref{222} gives a new identity:
\begin{exam}
\begin{align*}
&\begin{ytableau}
  2&2&2\\
  2\\
\end{ytableau}+\begin{ytableau}
  2&2\\
  2&2
\end{ytableau}\\
&=\sum_{\substack{(m_{ij})\in SSYT(\lambda)\\ (m)\in SSYT(e^1)}}\left[\frac{1}{m^2(m+m_{11})^2(m+m_{21})^2m_{12}^2}+\frac{2}{m(m+m_{11})^3(m+m_{21})^2m_{12}^2}\right.\\
&+\left.\frac{2}{m(m+m_{11})^2(m+m_{21})^3m_{12}^2}\right]+\begin{ytableau}
  2&2\\
  2\\
  2
\end{ytableau}+2\ \begin{ytableau}
  1&2\\
  3\\
  2
\end{ytableau}+2\ \begin{ytableau}
  1&2\\
  2\\
  3
\end{ytableau}+4\ \begin{ytableau}
  2&2\\
  1\\
  3
\end{ytableau}.
\end{align*}
\end{exam}
A similar calculation to the proof of Example \ref{222} for the product
\[\begin{ytableau}
  2&\pmb k\\
  2\\
\end{ytableau}\cdot\begin{ytableau}
  2
\end{ytableau}\]
with $\begin{ytableau}
  \pmb k
\end{ytableau}=\begin{ytableau}
  k_2&\cdots& k_r
\end{ytableau}$, where $k_r>1$, leads to the following:
\begin{exam}
\label{222k}
\begin{align*}
&\begin{ytableau}
  2&\pmb k\\
  2\\
\end{ytableau}\cdot\begin{ytableau}
  2
\end{ytableau}\\
&=\sum_{\substack{(m_{ij})\in SSYT(\lambda)\\ (m)\in SSYT(e^1)}}\prod_{j=2}^r\frac{1}{m_{1j}^{k_{j}}}\left[\frac{1}{m^2(m+m_{11})^2(m+m_{21})^2}+\frac{2}{m(m+m_{11})^3(m+m_{21})^2}\right.\\
&+\left.\frac{2}{m(m+m_{11})^2(m+m_{21})^3}\right]+2\ \begin{ytableau}
  2&\pmb k\\
  2\\
  2
\end{ytableau}+2\ \begin{ytableau}
  1&\pmb k\\
  3\\
  2
\end{ytableau}+2\ \begin{ytableau}
  1&\pmb k\\
  2\\
  3
\end{ytableau}+4\ \begin{ytableau}
  2&\pmb k\\
  1\\
  3
\end{ytableau}.
\end{align*}
\end{exam}
\section{Elementary factorial Schur multiple zeta functions}
\label{elementarysec}
As we observed in Section \ref{preex}, non-Schur multiple zeta values appear in the products of a Schur multiple zeta values of hook type and a multiple zeta values of Euler-Zagier type.
In this section, we describe how these values are expressed in terms of the modified Schur multiple zeta values. First, as a generalization of the Schur multiple zeta functions, we define a modified Hurwitz-type function as follows.
\begin{defi}
Let $\lambda=(\lambda_1,\{1\}^{k-1})$.
For given sets ${ \pmb s}=(s_{ij})\in T(\lambda,\mathbb{C}), { \pmb t}=(t_{i})\in T(e^r,\mathbb{C})$ of variables, we define
\[
\zeta_{\lambda,r}(\pmb s,\pmb t)=\sum_{\substack{(m_{ij})\in SSYT(\lambda)\\ (m_i)\in SSYT(e^r)}}\prod_{i=1}^r\frac{1}{m_{i}^{t_i}}\prod_{j=2}^{\lambda_1}\frac{1}{m_{1j}^{s_{1j}}}\prod_{i=1}^k\frac{1}{(m_r+m_{i1})^{s_{i1}}}.
\]
We call these {\it elementary factorial Schur multiple zeta functions}. Also, their convergent special values at positive integers are called {\it elementary factorial Schur multiple zeta values} (EFSMZVs).
\end{defi}
\begin{remark}
\label{matsumoto}
In \cite{mn}, the authors introduced modified zeta functions of root systems of type $A$: For $r>0$ and $0\le d\le r$,
\begin{equation}
    \label{rootsystem}
    \zeta_{r,d}^\bullet(\underbar{{\bf s}},A_r)=\left(\sum_{m_1=0}^\infty\cdots\sum_{m_d=0}^\infty\right)'\left(\sum_{m_{d+1}=1}^\infty\cdots\sum_{m_r=1}^\infty\right)\prod_{1\le i<j\le r+1}(m_i+\cdots+m_{j-1})^{-s(i,j)}.
\end{equation}
We understand that the components of the vector $\underbar{{\bf s}}=(s{(i,j)})$ 
is arranged according to the length of the corresponding term,
that is,
\begin{align*}
\underbar{{\bf s}}=(s{(1,2)},s{(2,3)},\ldots,s{(r,r+1)},s{(1,3)},s{(2,4)}&,\ldots,s{(r-1,r+1)},\\
&s{(1,r)}, s{(2,r+1)},s{(1,r+1)}).
\end{align*}
When $\lambda=(\lambda_1,1^{k-1})$, $r=\lambda_1+k-1+\ell$ and $d=\lambda_1-1$, \[s_{ij}=\left\{\begin{array}{cc}
s(\lambda_1-j+1,\lambda_1)    & 1< j\le \lambda_1,  \\
s(\lambda_1,\lambda_1+\ell+i-1)    & 1\le i\le k, 
\end{array}\right.\]
\[t_{i}=
s(\lambda_1+\ell-i+1,\lambda_1+\ell),\]
and the other $s(i,j)=0$,
then \[\zeta_{r,d}^\bullet(\underbar{{\bf s}},A_r)=\zeta_{\lambda,\ell}(\pmb s,\pmb t).\]
Therefore, we can regard the elementary factorial Schur multiple zeta function as a special case of (\ref{rootsystem}).
\end{remark}
We rewrite the example obtained by Theorem \ref{shuffle} in terms of these functions.
\begin{exam}[Rewrite of Example \ref{222k}]
\label{re222}
\begin{align*}\begin{ytableau}
  2&\pmb k\\
  2\\
\end{ytableau}\cdot\begin{ytableau}
  2
\end{ytableau}&=\zeta_{\lambda,1}\left(\begin{ytableau}
  2&\pmb k\\
  2
\end{ytableau},\begin{ytableau}2\end{ytableau}\right)+2\zeta_{\lambda,1}\left(\begin{ytableau}
  3&\pmb k\\
  2
\end{ytableau},\begin{ytableau}1\end{ytableau}\right)+2\zeta_{\lambda,1}\left(\begin{ytableau}
  2&\pmb k\\
  3
\end{ytableau},\begin{ytableau}1\end{ytableau}\right)\\
&+2\ \begin{ytableau}
  2&\pmb k\\
  2\\
  2
\end{ytableau}+2\ \begin{ytableau}
  1&\pmb k\\
  3\\
  2
\end{ytableau}+2\ \begin{ytableau}
  1&\pmb k\\
  2\\
  3
\end{ytableau}+4\ \begin{ytableau}
  2&\pmb k\\
  1\\
  3
\end{ytableau}.
\end{align*}
Next, we rewrite Example \ref{12112}.
For  $\lambda=(r,1)$, $\mu=(r,1,1)$, and $\nu=(r,1,1,1)$, we obtain another expression of Example \ref{12112}:
\begin{align*}\ytableausetup{boxsize=normal,aligntableaux=center}
&\begin{ytableau}
  1&\pmb k\\
  2
\end{ytableau}\times\begin{ytableau}
  1\\
  1\\
  2
\end{ytableau}\\
&=\sum_{\substack{(m_{ij})\in SSYT(\lambda)\\ (m_i)\in SSYT(e^3)}}\prod_{j=2}^r\frac{1}{m_{1j}^{k_{j}}}\left[\frac{1}{m_1m_2m_3^2(m_3+m_{11})(m_3+m_{21})^2}\right.\\
&\left.+\frac{1}{m_1m_2m_3(m_3+m_{11})^2(m_3+m_{21})^2}+\frac{2}{m_1m_2m_3(m_3+m_{11})(m_3+m_{21})^3}\right]\\
&+\sum_{\substack{(m_{ij})\in SSYT(\mu)\\ (m_i)\in SSYT(e^2)}}\prod_{j=2}^r\frac{1}{m_{1j}^{k_{j}}}\left[\frac{2}{m_1m_2(m_2+m_{11})(m_2+m_{21})^2(m_2+m_{31})^2}\right.\\
&\left.+\frac{4}{m_1m_2(m_2+m_{11})(m_3+m_{11})(m_3+m_{21})^3}\right]\\
&+\sum_{\substack{(m_{ij})\in SSYT(\nu)\\ (m_i)\in SSYT(e^1)}}\prod_{j=2}^r\frac{1}{m_{1j}^{k_{j}}}\left[\frac{3}{m_1(m_1+m_{11})(m_1+m_{21})(m_1+m_{31})^2(m_1+m_{41})^2}\right.\\
&\left.+\frac{6}{m_1(m_1+m_{11})(m_1+m_{21})(m_1+m_{31})(m_1+m_{41})^3}\right.\\
&\left.+\frac{1}{m_1(m_1+m_{11})(m_1+m_{21})^2(m_1+m_{31})(m_1+m_{41})^2}\right]\\
&+\begin{ytableau}
  1&\pmb k\\
  2\\
  1\\
  1\\
  2
\end{ytableau}+2\ \begin{ytableau}
  1&\pmb k\\
  1\\
  2\\
  1\\
  2
\end{ytableau}+4\ \begin{ytableau}
  1&\pmb k\\
  1\\
  1\\
  2\\
  2
\end{ytableau}+8\ \begin{ytableau}
  1&\pmb k\\
  1\\
  1\\
  1\\
  3
\end{ytableau}.
\end{align*} 
\end{exam}
\begin{exam}[Rewrite of Example \ref{12112}]
\begin{align*}\ytableausetup{boxsize=normal,aligntableaux=center}
\begin{ytableau}
  1&\pmb k\\
  2
\end{ytableau}\times\begin{ytableau}
  1\\
  1\\
  2
\end{ytableau}&=\zeta_{\lambda,3}\left(\begin{ytableau}
  1&\pmb k\\
  2
\end{ytableau},\begin{ytableau}1\\1\\2\end{ytableau}\right)+\zeta_{\lambda,3}\left(\begin{ytableau}
  2&\pmb k\\
  2
\end{ytableau},\begin{ytableau}1\\1\\1\end{ytableau}\right)+2\zeta_{\lambda,3}\left(\begin{ytableau}
  1&\pmb k\\
  3
\end{ytableau},\begin{ytableau}1\\1\\1\end{ytableau}\right)\\
&+2\zeta_{\mu,2}\left(\begin{ytableau}  
1&\pmb k\\
  2\\
  2
\end{ytableau},\begin{ytableau}1\\1\end{ytableau}\right)+4\zeta_{\mu,2}\left(\begin{ytableau}
  1&\pmb k\\
  1\\
  3
\end{ytableau},\begin{ytableau}1\\1\end{ytableau}\right)+3\zeta_{\nu,1}\left(\begin{ytableau}
  1&\pmb k\\
  1\\
  2\\
  2
\end{ytableau},\begin{ytableau}1\end{ytableau}\right)\\
&+6\zeta_{\nu,1}\left(\begin{ytableau}
  1&\pmb k\\
  1\\
  1\\
  3
\end{ytableau},\begin{ytableau}1\end{ytableau}\right)+\zeta_{\nu,1}\left(\begin{ytableau}
  1&\pmb k\\
  2\\
  1\\
  2
\end{ytableau},\begin{ytableau}1\end{ytableau}\right)\\
&+\begin{ytableau}
  1&\pmb k\\
  2\\
  1\\
  1\\
  2
\end{ytableau}+2\ \begin{ytableau}
  1&\pmb k\\
  1\\
  2\\
  1\\
  2
\end{ytableau}+4\ \begin{ytableau}
  1&\pmb k\\
  1\\
  1\\
  2\\
  2
\end{ytableau}+8\ \begin{ytableau}
  1&\pmb k\\
  1\\
  1\\
  1\\
  3
\end{ytableau}.
\end{align*} 
\end{exam}

In the following, we consider a relation  between the elementary factorial Schur multiple zeta functions and other multiple zeta functions.
Let $\mathcal Z_{ef}$ and $\mathcal Z_{sh}$ be the $\mathbb{Q}$-vector spaces generated by the elementary factorial  Schur multiple zeta values and by the Schur multiple zeta values, respectively.
Since we can decompose $\zeta_{\lambda,r}(\pmb s,\pmb t)$ into a sum of Schur multiple zeta functions, we find that $\mathcal Z_{ef}=\mathcal{Z}_{sh}$. Moreover, it is also known that $\mathcal Z_{sh}$ coincides the $\mathbb{Q}$-vector space $\mathcal Z$ generated by the multiple zeta values of Euler-Zagier type. Therefore, we have $\mathcal{Z}_{ef}=\mathcal{Z}$.
In fact, we can see the following examples.
\begin{exam}
\begin{align*}
\zeta_{\lambda,1}\left(\begin{ytableau}
  1&2\\
  2
\end{ytableau},\begin{ytableau}3\end{ytableau}\right)=&\begin{ytableau}3\end{ytableau}\cdot \begin{ytableau}1&2\\2\end{ytableau}-2\ \begin{ytableau}2\end{ytableau}\cdot \begin{ytableau}1&2\\3\end{ytableau}+3\ \begin{ytableau}\none&1&2\\1&4\end{ytableau}-\begin{ytableau}1&2\\3\\2\end{ytableau}\\
&-\begin{ytableau}2\end{ytableau}\cdot \begin{ytableau}2&2\\2\end{ytableau}+2\ \begin{ytableau}\none&2&2\\1&3\end{ytableau}-\begin{ytableau}2&2\\2\\2\end{ytableau}+\begin{ytableau}\none&3&2\\1&2\end{ytableau}-2\ \begin{ytableau}3&2\\1\\2\end{ytableau},\\
\zeta_{\lambda,2}\left(\begin{ytableau}
  1&2\\
  2
\end{ytableau},\begin{ytableau}1\\1\end{ytableau}\right)=&-\begin{ytableau}2\end{ytableau}\cdot \begin{ytableau}1&2\\2\end{ytableau}
+ \begin{ytableau}\none&\none&1&2\\1&1&2\end{ytableau}
-2\ \begin{ytableau}\none&1&2\\\none&1\\1&2\end{ytableau}
+3\ \begin{ytableau}1&2\\1\\1\\2\end{ytableau}.
\end{align*}
\end{exam}
\begin{proof}
The assertions of Example 4.4 come from direct calculation. Here, we only explain the details of the second assertion.
\begin{align*}
&\zeta_{\lambda,2}\left(\begin{ytableau}
  1&2\\
  2
\end{ytableau},\begin{ytableau}1\\1\end{ytableau}\right)\\
&=\sum_{\substack{(m_{ij})\in SSYT(\lambda)\\ (m_i)\in SSYT(e^2)}}\frac{1}{m_{1}m_2(m_2+m_{11})m_{12}^{2}(m_2+m_{21})^{2}}\\
&=\sum_{\substack{(m_{ij})\in SSYT(\lambda)\\ (m_i)\in SSYT(e^2)}}\frac{1}{m_{1}m_2m_{11}m_{12}^{2}(m_2+m_{21})^{2}}-\sum_{\substack{(m_{ij})\in SSYT(\lambda)\\ (m_i)\in SSYT(e^2)}}\frac{1}{m_{1}m_{11}(m_2+m_{11})m_{12}^{2}(m_2+m_{21})^2}\\
&=\sum_{\substack{(m_{ij})\in SSYT(\lambda)\\ (m_i)\in SSYT(e^2)}}\frac{1}{m_{1}m_2m_{11}m_{12}^{2}m_{21}(m_2+m_{21})}
-\sum_{\substack{(m_{ij})\in SSYT(\lambda)\\ (m_i)\in SSYT(e^2)}}\frac{1}{m_{1}m_{11}m_{12}^{2}m_{21}(m_2+m_{21})^2}\\
&-\sum_{\substack{(m_{ij})\in SSYT(\lambda)\\ (m_i),(m)\in SSYT(e^1)}}\frac{1}{m_{1}m_{11}(m_1+m+m_{11})m_{12}^{2}(m_1+m+m_{21})^2}\\
&=\sum_{\substack{(m_{ij})\in SSYT(\lambda)\\ (m_i)\in SSYT(e^2)}}\left[\frac{1}{m_{1}m_2m_{11}m_{12}^{2}m_{21}^2}
-\frac{1}{m_{1}m_{11}m_{12}^{2}m_{21}^2(m_2+m_{21})}\right]\\
&-\sum_{\substack{(m_{ij})\in SSYT(\lambda)\\ (m),(m_i)\in SSYT(e^1)}}\frac{1}{m_{1}m_{11}m_{12}^{2}m_{21}(m_1+m+m_{21})^2}\\
&-\sum_{\substack{(m_{ij})\in SSYT(\lambda)\\ (m_i),(m)\in SSYT(e^1)}}\frac{1}{m_{1}m_{11}(m+m_{11})m_{12}^{2}(m_1+m+m_{21})^2}\\
&+\sum_{\substack{(m_{ij})\in SSYT(\lambda)\\ (m_i),(m)\in SSYT(e^1)}}\frac{1}{m_{11}(m+m_{11})(m_1+m+m_{11})m_{12}^{2}(m_1+m+m_{21})^2}\\
&=\sum_{\substack{(m_{ij})\in SSYT(\lambda)\\ (m_i),(m)\in SSYT(e^1)}}\left[\frac{1}{m_{1}(m_1+m)m_{11}m_{12}^{2}m_{21}^2}
-\frac{1}{m_{1}m_{11}m_{12}^{2}m_{21}^2(m_1+m+m_{21})}\right]\\
&-\sum_{\substack{(m_{ij})\in SSYT(\lambda)\\ (m),(m_i)\in SSYT(e^1)}}\frac{1}{m_{1}m_{11}m_{12}^{2}m_{21}(m+m_{21})(m_1+m+m_{21})}\\
&+\sum_{\substack{(m_{ij})\in SSYT(\lambda)\\ (m),(m_i)\in SSYT(e^1)}}\frac{1}{m_{11}m_{12}^{2}m_{21}(m+m_{21})(m_1+m+m_{21})^2}\\
&-\sum_{\substack{(m_{ij})\in SSYT(\lambda)\\ (m_i),(m)\in SSYT(e^1)}}\frac{1}{m_{1}m_{11}(m+m_{11})m_{12}^{2}(m+m_{21})(m_1+m+m_{21})}\\
&+\sum_{\substack{(m_{ij})\in SSYT(\lambda)\\ (m_i),(m)\in SSYT(e^1)}}\frac{1}{m_{11}(m+m_{11})m_{12}^{2}(m+m_{21})(m_1+m+m_{21})^2}+\begin{ytableau}1&2\\1\\1\\2\end{ytableau}\\
&=\sum_{\substack{(m_{ij})\in SSYT(\lambda)\\ (m_i),(m)\in SSYT(e^1)}}\left[\frac{1}{m_{1}(m_1+m)m_{11}m_{12}^{2}m_{21}^2}
-\frac{1}{m_{1}m_{11}m_{12}^{2}m_{21}^2(m+m_{21})}\right.\\
&\left.+\frac{1}{m_{11}m_{12}^{2}m_{21}^2(m+m_{21})(m_1+m+m_{21})}\right]\\
&-\sum_{\substack{(m_{ij})\in SSYT(\lambda)\\ (m),(m_i)\in SSYT(e^1)}}\left[\frac{1}{m_{1}m_{11}m_{12}^{2}m_{21}(m+m_{21})^2}-\frac{1}{m_{11}m_{12}^{2}m_{21}(m+m_{21})^2(m_1+m+m_{21})}\right]\\
&+\begin{ytableau}1&2\\1\\1\\2\end{ytableau}-\sum_{\substack{(m_{ij})\in SSYT(\lambda)\\ (m_i),(m)\in SSYT(e^1)}}\frac{1}{m_{11}(m+m_{11})m_{12}^{2}(m+m_{21})^2}\left[\frac{1}{m_{1}}-\frac{1}{(m_1+m+m_{21})}\right]\\
&+\begin{ytableau}1&2\\1\\1\\2\end{ytableau}+\begin{ytableau}1&2\\1\\1\\2\end{ytableau}\\
&=\sum_{\substack{(m_{ij})\in SSYT(\lambda)\\ (m_i),(m)\in SSYT(e^1)}}\left[-\frac{1}{m_{1}^2m_{11}m_{12}^{2}m_{21}^2}+\frac{1}{m_{1}(m_1+m-1)m_{11}m_{12}^{2}m_{21}^2}\right.\\
&\left.-\frac{1}{m_{1}m_{11}m_{12}^{2}m_{21}^2(m+m_{21})}+\frac{1}{m_{11}m_{12}^{2}m_{21}^2(m+m_{21})(m_1+m+m_{21})}\right]\\
&-\begin{ytableau}\none&1&2\\\none&1\\1&2\end{ytableau}+3\ \begin{ytableau}1&2\\1\\1\\2\end{ytableau}
-\begin{ytableau}\none&1&2\\\none&1\\1&2\end{ytableau}\\
&=-\begin{ytableau}2\end{ytableau}\cdot \begin{ytableau}1&2\\2\end{ytableau}
+ \begin{ytableau}\none&\none&1&2\\1&1&2\end{ytableau}
-2\ \begin{ytableau}\none&1&2\\\none&1\\1&2\end{ytableau}
+3\ \begin{ytableau}1&2\\1\\1\\2\end{ytableau}.
\end{align*}
\end{proof}
Example \ref{222k} was given using integral expressions in Section \ref{2posets}. On the other hand, we can also give another proof using series expressions.
\begin{exam}[Another proof of Examples \ref{222k}, \ref{re222}]
We obtain the following three identities:
\begin{align*}
    \zeta_{\lambda,1}\left(\begin{ytableau}
  2&\pmb k\\
  2
\end{ytableau},\begin{ytableau}2\end{ytableau}\right)&=\begin{ytableau}
  2&\pmb k\\
  2\\
\end{ytableau}\cdot\begin{ytableau}
  2
\end{ytableau}-2\begin{ytableau}
  \none&2&\pmb k\\
  1&3
\end{ytableau}-2\begin{ytableau}
  \none&3&\pmb k\\
  1&2
\end{ytableau}+4\begin{ytableau}
  3&\pmb k\\
  1\\
  2
\end{ytableau}+2\begin{ytableau}
  2&\pmb k\\
  2\\
  2
\end{ytableau},\\
\zeta_{\lambda,1}\left(\begin{ytableau}
  3&\pmb k\\
  2
\end{ytableau},\begin{ytableau}1\end{ytableau}\right)&=\begin{ytableau}
  \none&3&\pmb k\\
  1&2
\end{ytableau}-2\begin{ytableau}
  3&\pmb k\\
  1\\
  2
\end{ytableau}-\begin{ytableau}
  2&\pmb k\\
  2\\
  2
\end{ytableau}-\begin{ytableau}
  1&\pmb k\\
  3\\
  2
\end{ytableau},\\
\zeta_{\lambda,1}\left(\begin{ytableau}
  2&\pmb k\\
  3
\end{ytableau},\begin{ytableau}1\end{ytableau}\right)&=\begin{ytableau}
  \none&2&\pmb k\\
  1&3
\end{ytableau}-\begin{ytableau}
  2&\pmb k\\
  2\\
  2
\end{ytableau}-2\begin{ytableau}
  2&\pmb k\\
  1\\
  3
\end{ytableau}-\begin{ytableau}
  1&\pmb k\\
  2\\
  3
\end{ytableau}.
\end{align*}
Substituting these three identities into the right-hand side of Example \ref{re222}, we complete the proof of Example \ref{222k}, as well as that of Example \ref{re222}.
\end{exam}
We can generalize $2$-labeled Schur posets to those for the corresponding elementary factorial Schur multiple zeta values, and define the shuffle product of these values as one generalization of Theorem \ref{shuffle}.
We note that the product can be decomposed into the sum of elementary factorial Schur multiple zeta values.
For example, we can write
\begin{align}
\label{shuffleexample}
    &\zeta_{\lambda,1}\left(\begin{ytableau}
  1&\pmb k\\
  2
\end{ytableau},\begin{ytableau}2\end{ytableau}\right)\times\begin{ytableau}
  2
\end{ytableau}\nonumber\\
&=2\zeta_{\lambda,2}\left(\begin{ytableau}
  1&\pmb k\\
  2
\end{ytableau},\begin{ytableau}2\\2\end{ytableau}\right)+4\zeta_{\lambda,2}\left(\begin{ytableau}
  1&\pmb k\\
  2
\end{ytableau},\begin{ytableau}1\\3\end{ytableau}\right)+2\zeta_{\lambda,2}\left(\begin{ytableau}
  2&\pmb k\\
  2
\end{ytableau},\begin{ytableau}1\\2\end{ytableau}\right)\nonumber\\
&\vspace{-6mm}\\
&+4\zeta_{\lambda,2}\left(\begin{ytableau}
  1&\pmb k\\
  3
\end{ytableau},\begin{ytableau}1\\2\end{ytableau}\right)+\zeta_{\lambda,2}\left(\begin{ytableau}
  2&\pmb k\\
  2
\end{ytableau},\begin{ytableau}2\\1\end{ytableau}\right)+2\zeta_{\lambda,2}\left(\begin{ytableau}
  1&\pmb k\\
  3
\end{ytableau},\begin{ytableau}2\\1\end{ytableau}\right)\nonumber\\
&+2\zeta_{\mu,1}\left(\begin{ytableau}
  1&\pmb k\\
  2\\
  2
\end{ytableau},\begin{ytableau}2\end{ytableau}\right)+4\zeta_{\mu,1}\left(\begin{ytableau}
  1&\pmb k\\
  1\\
  3
\end{ytableau},\begin{ytableau}2\end{ytableau}\right),\nonumber
\end{align}
where $\begin{ytableau}
  \pmb k
\end{ytableau}=\begin{ytableau}
  k_2&\cdots& k_r
\end{ytableau}$ and $\lambda=(r,1), \mu=(r,1,1)$.
Based on this property, we regard our shuffle product of Schur multiple zeta values (Theorem \ref{shuffle}) as that for more general values, those of elementary factorial type.
To formulate these, we redefine the $2$-labeled Schur posets.
\begin{defi}
\label{defi:sp2}
 Let $\begin{ytableau}
  \pmb k
\end{ytableau}=\begin{ytableau}
  k_2&\cdots& k_r
\end{ytableau}$.
\begin{enumerate}
\item A $2$-labeled Schur poset is a pair $((X,Y),\pmb k)=((X,Y),\pmb k;\delta_X)$ consisting of 
a poset $X$, a Young tableau $\pmb k$, and a labeling map $\delta\colon X\cup Y\to\{0,1\}$.
\item A $2$-labeled Schur poset $((X,Y),\pmb k)$ is called {admissible} 
if $\delta_Y(y)=1$ for all minimal $y\in Y$ and 
$\delta_X(x)=0$ for all maximal $x\in X$. 
\item For any posets $X$ and $Y$, we put 
\[\Delta(X,Y):=\bigl\{(t_\alpha)_{\alpha\in X\cup Y}\in[0,1]^{X\cup Y}\bigm| 
t_\alpha<t_\beta \text{ if } \alpha<\beta\bigr\}. \]
\item For admissible $2$-labeled Schur posets $I((X,Y),\pmb k)$, 
we define the associated integral by 
\begin{equation}\label{eq:I2}
I((X,Y),\pmb k):=\int_{\Delta(X,Y)}\prod_{\alpha\in X\cup Y}\omega_{\delta(\alpha)}(t_\alpha)\sum_{m_{1j}\in SSYT(h_{r-1})}\frac{1-t_{|X|}^{m_{12}}}{m_{12}^{k_{2}}\cdots m_{1r}^{k_r}},
\end{equation}
where $x$ is the minimal $x\in X$.
\end{enumerate}
\end{defi}
Also, we use Hasse diagrams to indicate $2$-labeled Schur posets, with vertices $\{\circ,\square\}$ and $\{\bullet,\blacksquare\}$ corresponding to $\delta(x)=0$ and $1$, respectively.
For example, a $2$-labeled Schur poset for $X=\{x_1,x_2,x_3,x_4\}$ and $Y=\{y_1,y_2,y_3\}$ with order $x_1>x_2>x_3>x_4>y_1>y_2>y_3$ and label $(\delta(x_1),\delta(x_2),\delta(x_3),\delta(x_4),$ $\delta(y_1),\delta(y_2),\delta(y_3))=(0,1,0,1,1,0,1)$ is represented as the diagram
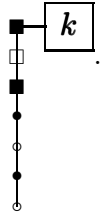
\[\begin{xy}
(7,4)*{\begin{ytableau}
  \pmb k
\end{ytableau}},
{(0,4) \ar @{-} (4,4)},
(0,4)*++{\scriptstyle \blacksquare},
(0,0)*++{\scriptstyle \square},
(0,-4)*++{\scriptstyle \blacksquare},
{(0,4) \ar @{-} (0,-8)},
{(0,-8) \ar @{{*}-o} (0,-12)},
{(0,-12) \ar @{-{*}} (0,-16)},
{(0,-16) \ar @{-o} (0,-20)},
\end{xy}.\]
This diagram corresponds to
\begin{align*}
&I\left(\begin{ytableau}
  2\\
  2
\end{ytableau},\begin{ytableau}2\\1\end{ytableau},\pmb k\right)\\
    &=\zeta_{\lambda,2}\left(\begin{ytableau}
  2&\pmb k\\
  2
\end{ytableau},\begin{ytableau}2\\1\end{ytableau}\right)\\
&=\int_{\Delta(X,Y)}\frac{dx_1}{x_1}\frac{dx_2}{1-x_2}\frac{dx_3}{x_3}\frac{dx_4}{1-x_4}\frac{dy_1}{1-y_1}\frac{dy_2}{y_2}\frac{dy_3}{1-y_3}\sum_{m_{1j}\in SSYT(h_{r-1})}\frac{1-x_4^{m_{12}}}{m_{12}^{k_{2}}\cdots m_{1r}^{k_r}},
\end{align*}
where $\Delta(X,Y)=\{1>x_1>x_2>x_3>x_4>y_1>y_2>y_3>0\}$.
We note that vertices of $X$ and $Y$ are distinguished. 
Considering the shuffle product of an elementary factorial Schur multiple zeta value and a multiple zeta value of Euler-Zagier type, it suffices to use the same vertices $\{\square,\blacksquare\}$ with $Y$ for multiple zeta values of Euler-Zagier type as before.
Then, one can generalize Theorem \ref{shuffle} as follows.
For $2$-labeled Schur posets $((X,Y),\pmb k)$ and $((Z,\emptyset),\emptyset)$, we define their direct sum $((X,Y)\amalg (Z,\emptyset),\pmb k)$: Its
underlying poset is the direct sum of finite sets $X\cup Y$ and $Z$ endowed with the partial order 
\[x< y \text{ in } (X,Y)\amalg (Z,\emptyset) \Leftrightarrow x, y\in X\cup Y \text{ and }x< y \text{ in } X\cup Y \text{ or
} x, y\in Z \text{ and }x< y \text{ in }Z.\]
\begin{thm}[Shuffle product of EFSMZV]
\label{shuffle2}
If $2$-labeled Schur posets $((X,Y),\pmb k)$ and $((Z,\emptyset),\emptyset)$ are admissible, then their direct sum $((X,Y)\amalg (Z,\emptyset),\pmb k)$ is also admissible and the following two properties hold.
\begin{enumerate}
    \item \[I((X,Y)\amalg (Z,\emptyset),\pmb k)=I((X,Y),\pmb k)I((Z,\emptyset),\emptyset),\]
    \item \[I((X,Y)\amalg (Z,\emptyset),\pmb k)\in \mathbb Z_{\ge0}[\mathcal{Z}_{ef}].\]
\end{enumerate}
\end{thm}
\begin{thm}[Explicit form of shuffle product of EFSMZV]
\label{shuffle3}
Let $\lambda=(\lambda_1,1^{r-1})$. For ${ \pmb k}=(k_{ij})\in T(\lambda,\mathbf{Z})$, ${ \pmb l}=(l_{i})\in T(e^s,\mathbf{Z})$ and ${ \pmb m}=(m_{i})\in T(e^t,\mathbf{Z})$, 
\begin{equation}
    \label{Shuffleprod}
    \zeta_{\lambda,s}\left(\pmb k,\pmb l\right)\cdot\zeta_{e^t}\left(\pmb m\right)=\sum_{p+q=t}\sum_{(\pmb u,\pmb v)}c_{\pmb u,\pmb v}\zeta_{\nu_p,q}\left(\pmb u,\pmb v\right),
\end{equation}
where $\nu_p=(\lambda_1,1^{r+s+p-1})$ and the sum $\displaystyle \sum_{(\pmb u,\pmb v)}$ means the summation over $({ \pmb u},{ \pmb v})=((u_{ij}), (v_{i}))$ such that for $j\ge2$, $u_{1j}=k_{1j}$ and $\zeta(v_1,\ldots,v_q,u_{11},\ldots,u_{(r+s+p)1})$'s appear in the summand of the shuffle product $\zeta(l_1,\ldots,l_s,k_{11},\ldots,k_{r1})\zeta(m_1,\ldots,m_t)$, and  $c_{\pmb u,\pmb v}$ are non-negative integers.
\end{thm}

Next, we calculate product (\ref{shuffleexample}) in a similar way to the original harmonic product. Then it holds that 
\begin{align*}
    \zeta_{\lambda,1}\left(\begin{ytableau}
  1&\pmb k\\
  2
\end{ytableau},\begin{ytableau}2\end{ytableau}\right)\times\begin{ytableau}
  2
\end{ytableau}=&\zeta_{\lambda,2}\left(\begin{ytableau}
  1&\pmb k\\
  2
\end{ytableau},\begin{ytableau}2\\2\end{ytableau}\right)+\zeta_{\lambda,1}\left(\begin{ytableau}
  1&\pmb k\\
  2
\end{ytableau},\begin{ytableau}4\end{ytableau}\right)\\
&+\sum_{\substack{(m_{ij})\in SSYT(\lambda)\\ (m_i)\in SSYT(e^2)}}\prod_{j=2}^r\frac{1}{m_{1j}^{k_{j}}}\frac{1}{m_1^2m_2^2(m_1+m_{11})(m_1+m_{21})^2}.
\end{align*}
Combining this identity with (\ref{shuffleexample}), we have
\begin{align*}
&\zeta_{\lambda,1}\left(\begin{ytableau}
  1&\pmb k\\
  2
\end{ytableau},\begin{ytableau}4\end{ytableau}\right)+\sum_{\substack{(m_{ij})\in SSYT(\lambda)\\ (m_i)\in SSYT(e^2)}}\prod_{j=2}^r\frac{1}{m_{1j}^{k_{j}}}\frac{1}{m_1^2m_2^2(m_1+m_{11})(m_1+m_{21})^2}\\
&=\zeta_{\lambda,2}\left(\begin{ytableau}
  1&\pmb k\\
  2
\end{ytableau},\begin{ytableau}2\\2\end{ytableau}\right)+4\zeta_{\lambda,2}\left(\begin{ytableau}
  1&\pmb k\\
  2
\end{ytableau},\begin{ytableau}1\\3\end{ytableau}\right)+2\zeta_{\lambda,2}\left(\begin{ytableau}
  2&\pmb k\\
  2
\end{ytableau},\begin{ytableau}1\\2\end{ytableau}\right)\\
&+4\zeta_{\lambda,2}\left(\begin{ytableau}
  1&\pmb k\\
  3
\end{ytableau},\begin{ytableau}1\\2\end{ytableau}\right)+\zeta_{\lambda,2}\left(\begin{ytableau}
  2&\pmb k\\
  2
\end{ytableau},\begin{ytableau}2\\1\end{ytableau}\right)+2\zeta_{\lambda,2}\left(\begin{ytableau}
  1&\pmb k\\
  3
\end{ytableau},\begin{ytableau}2\\1\end{ytableau}\right)\\
&+2\zeta_{\lambda_1,1}\left(\begin{ytableau}
  1&\pmb k\\
  2\\
  2
\end{ytableau},\begin{ytableau}2\end{ytableau}\right)+4\zeta_{\lambda_1,1}\left(\begin{ytableau}
  1&\pmb k\\
  1\\
  3
\end{ytableau},\begin{ytableau}2\end{ytableau}\right).
\end{align*}
\section{Further combinatorial or number theoretical problems}
In this section, we propose some further problems and briefly consider the corresponding examples.
\begin{prob}
Can we obtain an explicit formula for the coefficient $c_{\pmb u,\pmb v}$ in identity (\ref{Shuffleprod})?
\end{prob}
Euler \cite{euler} computed $c_{\pmb u,\pmb v}$ for the product of the Riemann zeta values $\zeta(k)\zeta(\ell)$, the most simple case, as 
\[\zeta(k)\zeta(\ell)=\sum_{i=1}^k\binom{k+\ell-i-1}{\ell-1}\zeta(k+\ell-i,i)+\sum_{i=1}^\ell\binom{k+\ell-i-1}{k-1}\zeta(k+\ell-i,i).\]
Recently, Guo and Xie obtained an explicit shuffle product formula in a very general setting and gave the coefficient $c_{\pmb u,\pmb v}$ for the product of $\zeta(k_1,k_2)\zeta(\ell)$ and $\zeta(k_1,k_2)\zeta(\ell_1,\ell_2)$ \cite{gx} explicitly. For example, for $k_1\ge1$ and $k_2,\ell\ge2$, we have
\begin{align*}
    \zeta(k_1,k_2)\zeta(\ell)=&\sum_{\substack{i_1\ge1,i_2\ge2\\i_1+i_2=k_2+\ell}}\binom{i_2-1}{\ell-1}\zeta(k_1,i_1,i_2)\\
    &+\sum_{\substack{i_1,i_2\ge1,i_3\ge2\\i_1+i_2+i_3=k_1+k_2+\ell}}\binom{i_3-1}{k_2-1}\left[\binom{i_2-1}{k_1-i_1}+\binom{i_2-1}{k_1-1}\right]\zeta(i_1,i_2,i_3).
\end{align*}
In \cite{elo,ew,lgm}, the authors studied the product of multiple zeta functions $\zeta(k_1,\ldots,k_r)$ with strings of $1$'s. More preciously, Eie and Wei computed the coefficient of the product $\zeta(\{1\}^j,m+r-\ell-2)\zeta(\{1\}^{k-j},q-m+\ell+2)$ \cite{ew}. Eie, Liaw, and Ong generalized this result to the product of multiple zeta values $\zeta(\{1\}^a,b+2)$ with non-negative integers $a,b$ \cite{elo}. As another generalization of \cite{ew}, Lie, Guo, and Ma also calculated the coefficient of the product $\zeta(\{1\}^i,\ell+1)\zeta(\{1\}^{j},m+1,\{1\}^k,n+1)$ with two strings of 1's \cite{lgm}.
Further, in \cite{lq}, Li and Qin expressed the product $\zeta(k_1+1,\ldots,k_r+1)\zeta(\ell_1+1,\ldots,\ell_s+1)$ with positive integers $k_i,\ell_i$ explicitly.

\begin{prob}
Can we obtain an explicit form of the decomposition of an elementary factorial Schur Multiple zeta function into a sum of (skew) Schur multiple zeta functions?
\end{prob}
\begin{exam}Let $\lambda=(\lambda_1,1^{r-1})$ and $\pmb k\in T(\lambda,\mathbb Z)$.
\begin{align*}
  \zeta_{\lambda,1}\left(
  \pmb k,\begin{ytableau}1\end{ytableau}\right)&=\begin{ytableau}
  \none&k_{11}&\cdots&k_{1\lambda_1}\\
\none&\vdots\\
  1&k_{1r}
\end{ytableau}-\sum_{i=1}^{r}\sum_{\substack{\ell_i+\ell_{i+1}=1+k_{i1}\\
\ell_i,\ell_{i+1}\ge1,\ell_{r+1}\ge2\\
\ell_1=k_{11},\ldots,\ell_{i-1}=k_{(i-1)1}\\
\ell_{i+2}=k_{(i+1)1},\ldots,\ell_{r+1}=k_{r1}\\}}\begin{ytableau}
  \ell_1&k_{12}&\cdots&k_{1\lambda_1}\\
  \vdots\\
  \ell_{r+1}
\end{ytableau}
\end{align*}
\end{exam}

Let $\pmb t_j=(t_{j1},\ldots,t_{jr_j})\in T(e^{r_j},\mathbb C)$, $\pmb s=(s_{ij})\in T(\lambda,\mathbb C)$. 
One may define the ``$n$-tuple'' elementary factorial Schur multiple zeta functions by
\begin{align*}
    \zeta(\pmb s,\pmb t_1,\pmb t_{2},\ldots,\pmb t_n)=&\sum_{\substack{(m_{ij})\in SSYT(\lambda)\\ (\mu_{1i})\in SSYT(e^{r_1})\\
\vdots\\
\\ (\mu_{ni})\in SSYT(e^{r_n})}}\left[\prod_{a=1}^n\prod_{b=1}^{r_a}\frac{1}{(\mu_{nr_n}+\cdots+\mu_{(a+1)r_{a+1}}+\mu_{ab})^{t_{ab}}}\prod_{c=2}^{\lambda_1}\frac{1}{m_{1c}^{s_{1c}}}\right.\\
&\times \left.\prod_{d=1}^k\frac{1}{(\mu_{nr_n}+\cdots+\mu_{1r_1}+m_{i1})^{s_{i1}}}\right].
\end{align*}
\begin{prob}
Do the ``$n$-tuple'' elementary factorial Schur multiple zeta values/functions have any interesting properties?
\end{prob}

\begin{remark}
If we denote by $\mathcal Z_{eh}^n$ the $\mathbb{Q}$-vector space generated by the ``$n$-tuple'' elementary factorial Schur multiple zeta values, then we find $\mathcal Z_{eh}^n=\cdots=\mathcal Z_{eh}^1=\mathcal Z$.
\end{remark}

\section*{Acknowledgement}
The authors deeply express their sincere gratitude to Professor Kohji Matsumoto for pointing out Remark \ref{matsumoto}.
The first author was supported by a Grant-in-Aid for Scientific Research (C) (Grant Number: JP18K03223).


\begin{thebibliography}{99}
\bibitem[ELO]{elo} M. Eie, W.-C. Liaw, Y. L. Ong. The decomposition theorem of products of multiple zeta values of height one
Int. J. Number Theory, {\bf 12} (2016), 15--25.
\bibitem[EM]{em} S. Eilenberg and S. Mac Lane. On the Groups $H(\Pi,n)$, I. Annals of Mathematics,
{\bf 58} (1953), 55-106.
\bibitem[EW]{ew} M. Eie and C.-S. Wei.
Generalizations of Euler decomposition and their applications.
J. Number Theory, {\bf 133} (2013), 2475--2495.
\bibitem[E]{euler} L. Euler. Meditationes circa singulare serierum genus. Novi. Comm. Acad. Sci. Petropolitanae, {\bf 20} (1775), 140--186.
\bibitem[GX]{gx} L. Guo and B. Xie.
Explicit double shuffle relations and a generalization of Euler's decomposition formula. J. Algebra, {\bf 380} (2013), 46--77.
\bibitem[H]{hoffman} M. Hoffman. The algebra of multiple harmonic series. J. Algebra, {\bf 194} (2) (1997), 477--495.
\bibitem[HO]{hoffmanohno} M. Hoffman and Y. Ohno. Relations of multiple zeta values and their algebraic expression. J. Algebra, {\bf 262} (2003), 332--347.
\bibitem[LGM]{lgm} P. Lei and L. Guo and B. Ma. Applications of shuffle product to restricted decomposition formulas for multiple zeta values. J. Number Theory, {\bf 144} (2014), 219--233.
\bibitem[LQ]{lq} Z. Li and C. Qin. Shuffle product formulas of multiple zeta values. J. Number Theory, {\bf 171} (2017), 79--111.
\bibitem[MN]{mn} K. Matsumoto and M. Nakasuji. Expressions of Schur multiple zeta-functions of anti-hook type by zeta-functions of root systems. Publ. Math. Debrecen {\bf 98} (2021), no. 3-4, 345--377.
   \bibitem[NN]{nn} N. Nakamura and M. Nakasuji. Schur type poly-Bernoulli numbers. preprint
  \bibitem[NPY]{npy} M. Nakasuji, O. Phuksuwan and Y. Yamasaki. On Schur multiple zeta functions: A combinatoric generalization of multiple zeta functions, Advances in Mathematics, {\bf 333} (2018), 570--619. 
  \bibitem[{N}{T}]{nt}M. Nakasuji and W. Takeda,
The Pieri formulas for hook type Schur multiple zeta functions,
arXiv : 2105.12418.
\bibitem[Y]{y} S. Yamamoto, Multiple zeta-star values and multiple integrals. RIMS
Kokyuroku Bessatsu B68 (2017), 3--14.
\end{thebibliography}
\end{document}